\documentclass{gtmon_a}
\pdfoutput=1

\usepackage{amscd}

\proceedingstitle{Proceedings of the Nishida Fest (Kinosaki 2003)}
\conferencestart{28 July 2003}
\conferenceend{8 August 2003}
\conferencename{International Conference in Homotopy Theory}
\conferencelocation{Kinosaki, Japan}

\editor{Matthew Ando}
\givenname{Matthew}
\surname{Ando}

\editor{Norihiko Minami}
\givenname{Norihiko}
\surname{Minami}

\editor{Jack Morava}
\givenname{Jack}
\surname{Morava}

\editor{W Stephen Wilson}
\givenname{W Stephen}
\surname{Wilson}

\title{The stable braid group and the determinant of the Burau representation}

\author{F\,R Cohen}
\givenname{F\,R}
\surname{Cohen}
\address{Department of Mathematics\\
University of Rochester\\\newline
Rochester NY 14627\\
USA}
\email{cohf@math.rochester.edu}
\urladdr{}

\author{J Pakianathan}
\givenname{J}
\surname{Pakianathan}
\email{jonpak@math.rochester.edu}
\urladdr{}

\volumenumber{10}
\issuenumber{}
\publicationyear{2007}
\papernumber{5}
\startpage{117}
\endpage{129}

\doi{}
\MR{}
\Zbl{}

\arxivreference{math.AT/0509577}


\received{15 January 2004}
\revised{5 June 2005}
\accepted{}
\published{29 January 2007}
\publishedonline{29 January 2007}
\proposed{}
\seconded{}
\corresponding{}
\version{}

\makeatletter
\def\cnewtheorem#1[#2]#3{\newtheorem{#1}{#3}[section]
\expandafter\let\csname c@#1\endcsname\c@thm}

\def\<{\langle}
\def\>{\rangle}

\let\xysavmatrix\xymatrix
\def\xymatrix{\disablesubscriptcorrection\xysavmatrix}
\AtBeginDocument{\let\bar\wbar\let\tilde\wtilde\let\hat\what}

\newtheorem{thm}{Theorem}[section]
\cnewtheorem{cor}[thm]{Corollary}
\cnewtheorem{lem}[thm]{Lemma}
\cnewtheorem{prop}[thm]{Proposition}

\theoremstyle{definition}
\cnewtheorem{defin}[thm]{Definition}
\cnewtheorem{exm}[thm]{Example}

\theoremstyle{remark}
\newtheorem*{rem}{Remark}

\makeatother  %  move after \newtheorem block

\makeop{rep}
\makeop{Rep}
\makeop{Hom}
\makeop{Irr}
\makeop{Top}
\makeop{Sym}

\begin{document}

\begin{htmlabstract}
<p class="noindent">
This article gives certain fibre bundles associated to the braid
groups which are obtained from a translation as well as conjugation
on the complex plane. The local coefficient systems on the level of
homology for these bundles are given in terms of the determinant of
the Burau representation.</p>

<p class="noindent">
De Concini, Procesi, and Salvetti [Topology 40 (2001) 739&ndash;751]
considered the cohomology of the n<sup>th</sup>
braid group B<sub>n</sub> with local coefficients
obtained from the determinant of the Burau representation,
H<sup>*</sup>(B<sub>n</sub>;<b>Q</b>[t<sup>&plusmn;1</sup>]). They show
that these cohomology groups are given in terms of cyclotomic fields.</p>

<p class="noindent">
This article gives the homology of the stable
braid group with local coefficients obtained from the
determinant of the Burau representation. The main result is an isomorphism
H<sub>*</sub>(B<sub>&infin;</sub>;<b>F</b>[t<sup>&plusmn;1</sup>])&rarr;H<sub>*</sub>(&Omega;<sup>2</sup>S<sup>3</sup>&lang;3&rang;;<b>F</b>)
for any field <b>F</b> where &Omega;<sup>2</sup>S<sup>3</sup>&lang;3&rang;
denotes the double loop space of the 3&ndash;connected cover of
the 3&ndash;sphere. The methods are to translate the structure of
H<sub>*</sub>(B<sub>n</sub>;<b>F</b>[t<sup>&plusmn;1</sup>]) to one
concerning the structure of the homology of certain function spaces
where the answer is computed.</p>
\end{htmlabstract}

\begin{webabstract}
This article gives certain fibre bundles associated to the braid
groups which are obtained from a translation as well as conjugation
on the complex plane. The local coefficient systems on the level of
homology for these bundles are given in terms of the determinant of
the Burau representation.

De Concini, Procesi, and Salvetti considered the
cohomology of the $n$th braid group $B_n$ with local coefficients
obtained from the determinant of the Burau representation,
$H^*(B_n;\mathbb{Q}[t^{\pm 1}])$. They show that these cohomology
groups are given in terms of cyclotomic fields.

This article gives the homology of the stable braid group
with local coefficients obtained from the determinant
of the Burau representation. The main result is an
isomorphism
$$H_*(B_{\infty};\mathbb{F}[t^{\pm 1}])\to
  H_*(\Omega^2S^3\langle3\rangle;\mathbb{F})$$
for any field $\mathbb{F}$ where
$\Omega^2S^3\langle3\rangle$ denotes the double loop space of the
$3$--connected cover of the $3$--sphere. The methods are to translate
the structure of $H_*(B_n;\mathbb{F}[t^{\pm 1}])$ to one concerning the
structure of the homology of certain function spaces where the answer
is computed.
\end{webabstract}

\begin{abstract}
This article gives certain fibre bundles associated to the braid
groups which are obtained from a translation as well as conjugation
on the complex plane. The local coefficient systems on the level of
homology for these bundles are given in terms of the determinant of
the Burau representation.

De Concini, Procesi, and Salvetti \cite{DPS} considered the
cohomology of the $n$th braid group $B_n$ with local coefficients
obtained from the determinant of the Burau representation,
$H^*(B_n;\mathbb{Q}[t^{\pm 1}])$. They show that these cohomology
groups are given in terms of cyclotomic fields.

This article gives the homology of the stable braid group
with local coefficients obtained from the determinant
of the Burau representation. The main result is an
isomorphism
$$H_*(B_{\infty};\mathbb{F}[t^{\pm 1}])\to
  H_*(\Omega^2S^3\langle3\rangle;\mathbb{F})$$
for any field $\mathbb{F}$ where
$\Omega^2S^3\langle3\rangle$ denotes the double loop space of the
$3$--connected cover of the $3$--sphere. The methods are to translate
the structure of $H_*(B_n;\mathbb{F}[t^{\pm 1}])$ to one concerning the
structure of the homology of certain function spaces where the answer
is computed.
\end{abstract}

\maketitle

\section{Introduction}
\label{sec1}

This article gives certain fibre bundles associated to the braid
groups which are obtained from a natural action of translation, and
conjugation on the complex plane. The local coefficient systems on
the level of homology for these bundles are given in terms of the
determinant of the Burau representation. Furthermore, the following
are addressed:

\begin{enumerate}
\item The total space of each of the fibre bundles is a
$K(\pi,1)$ with local coefficients in homology given by the
determinant of the Burau representation. These are related to work
of C De Concini, C Procesi, and M Salvetti \cite{DPS} who analyze
the cohomology of the braid groups with local coefficients given by
the determinant of the Burau representation. These bundles arise by
pulling back a bundle given below which is obtained from an action
of the integers on the complex plane $\mathbb C$ with the integral
points $(n,0)$ deleted.

    \item The homology of the bundles obtained by the standard stablization of the
    braid groups is given by the homology for the $3$--connected cover of the double loop space
of the $3$--sphere. Thus homology is considered here rather than
cohomology as in the case of \cite{DPS}.

    \item The last section of this article is a short synopsis of a recurring
    relationship concerning the (co)homology of the braid groups with certain
    related local coefficient systems. Here the resulting groups are
    given by cyclotomic fields, as well
    as classical modular forms (see Furusawa--Tezuka--Yagita \cite{FTY}
    and Cohen \cite{C2}). Several related problems are posed.
\end{enumerate}

Let $\mathbb F$ denote a field. The ring of Laurent polynomials
$$\mathbb Z[t^{\pm 1}] \otimes_{\mathbb Z} \mathbb F = \mathbb F[t^{\pm 1}]$$
is given by the elements
$$ \Sigma_{ -L \leq i \leq L} a_{i}t^{i} = a_{-L}t^{-L} +
a_{-L+1}t^{-L+1} + \cdots + a_0 + a_{1}t^{1} + \cdots+ a_{L}t^{L}$$
for all $ 0 \leq L < \infty $.

Artin's presentation for  $n$th braid group $B_n$ (see Birman \cite{B}) is given
by generators  $\sigma_j$ with $1 \leq i \leq n-1$ with relations
\begin{itemize}
    \item $\sigma_i \sigma_j = \sigma_j \sigma_i$ if $|i-j|\geq
    2$, and
    \item $\sigma_i \sigma_{i+1}\sigma_i = \sigma_{i+1}\sigma_i \sigma_{i+1}.$
\end{itemize}
Observe that the braid group acts on the complex numbers $\mathbb C$
by the formula
$$\sigma_i(z) = 1 + \chi(z)$$
where $\chi(z)$ denotes
the complex conjugate of $z$. This action restricts to one on
\begin{enumerate}
  \item the space $X$ which is the complement of the standard
lattice $\mathbb Z + i \mathbb Z$, the Gaussian integers in $\mathbb
C$, and
  \item the space $Y$ which is
the complement of the integral points $\{(\lambda,0)| \lambda \in
\mathbb Z\}$ in $\mathbb C$.
\end{enumerate}
In addition, the spaces $X$, and $Y$ are both
$K(\pi,1)$'s for which $\pi$ is a countably infinitely generated
free group. Furthermore a choice of generators of $\pi_1(Y)$ is
naturally indexed by powers of the indeterminate $t$ in $\mathbb
Z[t^{\pm 1}]$ with the abelianization of $\pi_1(Y)$ isomorphic to
the underlying structure of $ \mathbb Z[t^{\pm 1}]$ as an abelian
group. One choice of generators of $\pi_1(X)$ is naturally indexed
by the points in the standard integral lattice for the complex
numbers.

The braid group acts on $\mathbb Z[t^{\pm 1}]$ where each $\sigma_i$
acts via multiplication by $-t$. This action is obtained from the
determinant of the classical Burau representation of $\sigma_i$,
$$b\co B_n \to\ GL(n, \mathbb Z[t^{\pm 1}])$$
with $b(\sigma_i)$ given by the matrix
$$\begin{pmatrix}
  1 & 0  & . & . & . & 0 \\
  0 & 1 & 0 & 0 & 0 & 0 \\
  . & 0 & a_{i,i} & a_{i,i+1} & 0 & . \\
  . & 0 & a_{i+1,i} & a_{i+1,i+1} & 0 & . \\
  . & . & 0 & 0 & 1 & . \\
  0 & 0 & . & . & 0 & 1
\end{pmatrix}$$
where $a_{i,i} = 0$, $a_{i+1,i} = 1$, $a_{i,i+1} = t$, and
$a_{i+1,i+1} = 1-t$. Thus $\mathbb Z[t^{\pm 1}]$ is a $B_n$--module
where $\sigma_i$ acts by multiplication by the determinant $d$ of
$b(\sigma_i)$ with
$$d(b(\sigma_i)) = -t.$$
This module structure is
precisely that given by De Concini, Procesi and Salvetti \cite{DPS}
after tensoring with the rational numbers.

The main result of \cite{DPS} is a computation of the cohomology of
the $n$th Artin braid group $B_n$ with coefficients in $\mathbb
Q[t^{\pm 1}]$ with the above action. As in \cite{DPS}, let
$\mu_d(q)$ denote the cyclotomic polynomial in $q$ having as roots
the primitive $d$th roots of unity. Let $R_q$ denote the
$B_{m+1}$--module given by the action on the ring
$$R_q =\mathbb Q[q^{\pm 1}]$$
of Laurent polynomials over the rational numbers in the variable $q$
which is defined by mapping each standard generator $\sigma_i$ to
the operator of multiplication  by $-q$. Define
$$\{d\}= \mathbb Q[q^{\pm 1}]/\mu_d(q)=\mathbb Q[q]/\mu_d(q),$$
the
cyclotomic field of $d$ roots of unity.

\begin{thm}[De Concini--Procesi--Salvetti \cite{DPS}]
\label{thm:procesi.deconcini.salvetti}
If $hi=m+1$ or $hi=m $ then
$$H^{i(h-2)+1}(B_{m+1},R_q)=\{h\}.$$
All other cohomology groups are trivial.
\end{thm}

The main result of this article addresses the homology of the stable
braid group $B_{\infty}$ with coefficients in  $\mathbb Z[t^{\pm
1}]$ rather than cohomology with analogous coefficients. The main
point here is that there are natural geometric bundle
interpretations of these results.

That is, consider the Borel construction or homotopy orbit space
$$EB_n \times_{B_n}Y,$$
the total space of a fibre bundle over $BB_n = K(B_n,1)$ with fibre
$Y = K(\pi_1(Y),1)$. Consequently, the space $EB_n \times_{B_n}Y$ is
a $K(\pi,1)$, and is equipped with natural ``stabilization" maps
$EB_n \times_{B_n}Y \to EB_{n+1} \times_{B_{n+1}}Y$ with colimit
denoted $EB_{\infty} \times_{B_{\infty} }Y$. The next proposition
records the local coefficient system for the bundle projection $EB_n
\times_{B_n}Y \to BB_n$.

\begin{prop}\label{prop:local.coefficients}
Assume that $ n \geq 2$. The action of $B_n$ on $H_1(Y)$ specified
by the local coefficient system in homology for the bundle
projection $EB_n \times_{B_n}Y \to BB_n$ is given by multiplication
by the determinant of the Burau representation: $\sigma_i(t^j) =
-t^{j+1}$.
\end{prop}

The homology of the space $EB_{\infty} \times_{B_{\infty} }Y$ with
coefficients in a field $\mathbb F$ is stated in the next theorem
for which the space $\Gamma(S^1)$ denotes $S^1 \vee S^1 = K(F_2,1)$
with $F_2$ a free group on two letters.

\begin{thm}\label{thm:bundles}

Assume that $ n \geq 2$.
\begin{enumerate}
    \item There is an induced map
$$\rho\co EB_{\infty} \times_{B_{\infty} }Y \to \Gamma(S^1)
  \times \Omega^2S^3\<3\>$$
which is a homology isomorphism
for any field $\mathbb F$ (with trivial action) where $\Gamma(S^1)$
is a $K(F_2,1)$. Thus the suspension of $\rho$, $\Sigma(\rho)$, is a
homotopy equivalence.

    \item There are isomorphisms
\begin{eqnarray*}
H_*(B_{\infty}; \mathbb F) \oplus H_*(B_{\infty}; \mathbb F[t^{\pm 1}])
&\to& H_*(EB_{\infty} \times_{B_{\infty} }Y; \mathbb F),\\
H_*(B_{\infty}; \mathbb F) \oplus H_*(B_{\infty}; \mathbb F[t^{\pm 1}])
&\to& H_*( \Gamma(S^1); \mathbb F) \otimes H_*( \Omega^2 S^3 \<3\>; \mathbb
F),\\
H_*(B_{\infty};\mathbb F[t^{\pm 1}]) &\to& H_{*}(\Omega^2 S^3 \<3\>;\mathbb F).
\end{eqnarray*}

\item If $i \geq 1$, there are isomorphisms $ H_{i+1}(EB_n
\times_{B_n}Y; \mathbb Q) \to\  H_i(B_n; \mathbb Q[t^{\pm 1}])$ for
$i > 0$.
\end{enumerate}
\end{thm}

This theorem is proven by analyzing properties of  fibre bundles
described $EB_n \times_{B_n}Y$ in \fullref{sec2}. In addition, the integer
homology of $\Omega^2 S^3\langle 3 \rangle$ is well-understood.  For
example, the reduced homology is entirely torsion for which the
$p$--torsion is of order exactly $p$.

The mod--$2$ homology is a polynomial ring with generators $x_i$ of
degree $i$ where $i$ is either 2 or $2^{n-1}$ for all $n>1$.  The first
Bockstein $\beta$ satisfies $\beta(x_3)=x_2$ and $\beta(x_{2^n-1})
= x^2_{2^{n-1}-1}$ if $n > 2$. If p is an odd prime, the mod-$p$
homology is a tensor product of (i) an exterior algebra with generators
$y_{2p^n-1}$ of degree $2p^n-1$, $n \geq 1$ and (ii) a polynomial algebra
with generators $x_{2p^n-2}$  of degree $2p^n-2$ for $n \geq 1$. The
first Bockstein  $\beta$ satisfies $\beta(y_{2p^n-1}) = x_{2p^n-2}$. One
reference is \cite{C}, but the structure as an algebra was first worked
out by J\,C Moore.

\begin{rem}
\label{remark:sections}
Filippo Callegaro has shown
that the bundles $EB_n \times_{B_n}Y \to BB_n$ admit cross-sections.
Thus if $i \geq 1$, there are isomorphisms
$$H_{i+1}(EB_n \times_{B_n}Y; \mathbb Z) \to H_{i+1}(B_n; \mathbb
Z)\oplus H_i(B_{n}; \mathbb Z[t^{\pm 1}]).$$
In addition, let
$\mathbb Z[t^{\pm 1}]^*$ denote the hom-dual of $\mathbb Z[t^{\pm
1}]$. The existence of a cross-section implies that there are
isomorphisms
$$H^{i+1}(EB_n \times_{B_n}Y; \mathbb Z) \to
H^{i+1}(B_n; \mathbb Z)\oplus H^i(B_{n}; \mathbb Z[t^{\pm 1}]^*).$$
Callegaro has also determined the precise homology of
$EB_n \times_{B_n}Y$ in \cite{callegaro}.
\end{rem}

The authors would like to congratulate Goro Nishida on this happy
occasion of his 60th birthday. His work has been very
interesting as well as very stimulating for the authors, especially
an early paper which greatly encouraged one of the authors. In
addition, the authors would like to thank Filippo Callegaro as well
as the referee for suggestions concerning this article.
,%MA Added this from \thanks
Cohen was partially supported by the National Science Foundation
under Grant No. 9704410 and CNRS-NSF Grant No. 17149.

%%SECTION 2
\section{Bundles associated to the determinant of the Burau
representation}
\label{sec2}

Regard $\mathbb Z[t^{\pm 1}]$ as a $B_n$--module obtained from the
determinant of the Burau representation as above. Let $EG$ denote
Milnor's contractible ``universal space" on which the discrete group
$G$ acts freely and properly discontinuously. The purpose of the
next lemma is to give elementary properties stated earlier.

\begin{lem}\label{lem:Action}
\begin{enumerate}
\item There is an action of $B_n$ on the complex line $\mathbb C$
given by
$$\sigma_j(z) = 1 +  \chi(z).$$

\item This action of $B_n$ on $\mathbb C$ restricts to an action
on the spaces $X$ and $Y$.

\item The associated  Borel constructions $EB_n \times_{B_n}W$,
where $W$ is either $X$, or $Y$ are $K(\pi,1)$'s. In case $W = Y$,
the group $\pi$ is given by a group extension
$$1 \to \pi_1(Y) \to \pi \to\ B_n \to\ 1$$
where $\pi_1(Y)$ is a free group with generators indexed by $t^q$
for $q$ in $\mathbb Z$ and where the action of $B_n$ on $H_1(Y)$ is
given by the action of the determinant of the Burau representation.
\end{enumerate}
\end{lem}

\begin{proof}
An action of $B_n$ on the complex numbers $\mathbb C$ is specified
by the composite
\[
\begin{CD}
B_n @>{A}>> \mathbb Z @>{\Theta}>>  \Top(\mathbb C)
\end{CD}
\] where $\Top(\mathbb C)$ is the homeomorphism group of $\mathbb
C$,
$$\Theta(1)(z) = 1 + \chi(z)$$
and
$$A\co B_n \to \mathbb Z$$
is given by abelianization with $A(\sigma_i) = 1$ for all $i$. Thus
$$\sigma_j(z) = 1 + \chi(z).$$
That this action restricts to an action on $X$, and $Y$ follows at
once as the standard integral lattice as well as the integral points
$\{(n,0)| n \in \mathbb Z\}$ in $\mathbb C$ are preserved by this
action. The second part follows.

Notice that there is a fibration
$$EB_n \times_{B_n} W  \to\  BB_n$$
with fibre $W$. Since $BB_n$, $X$ and $Y$
are  $K(\pi,1)'s$, so is the associated total space. Finally, the
group extension arises by applying the fundamental group functor to
the fibration $EB_n \times_{B_n} W  \to\  BB_n$ where all spaces are
$K(G,1)'s$. That the local coefficient system in homology for this
fibration is given by the action of the determinant of the Burau
representation is verified below in \fullref{lem:local.coefficients.again}.
The third part follows.
\end{proof}

Observe that the first homology group of $X$, $H_1(X)$, is a
countably infinite direct sum of copies of the integers indexed by
the points in the standard standard integral lattice in $\mathbb C$.
Similarly, the first homology group of $Y$, $H_1(Y)$, is a countably
infinite direct sum of copies of the integers indexed by the
integers. Furthermore, a basis for $H_1(Y)$ is given by $\iota_n$
for $n$ in $\mathbb Z$ for which $\iota_n$ is given by the image of
a fixed fundamental cycle of a small circle embedded in $Y$ with
center $(n,0)$, and radius $1/2$.

The next lemma directly implies \fullref{prop:local.coefficients} concerning the local coefficient
system in homology for the bundle projection $EB_n \times_{B_n}Y \to
BB_n$.

\begin{lem}\label{lem:local.coefficients.again}
There is an isomorphism of $B_{\infty}$--modules
$$\Theta\co \mathbb Z[t^{\pm 1}] \to\ H_1(Y)$$
given by
$\Theta(t^n) = \iota_n$ where $\iota_n$ is the homology class of a
small circle in $Y$ with center $(n,0)$, and radius $1/2$.
\end{lem}

\begin{proof}
That $\Theta$ is a homology isomorphism follows by inspection.
Notice that $\sigma_i$ sends $t^n$ to $(-t)t^n = -t^{n+1}$ in
$\mathbb Z[t^{\pm 1}]$. Furthermore, $\sigma_i$ sends the small
circle with center $(n,0)$, and radius $1/2$ homeomorphically to a
small circle with center $(n+1,0)$, and radius $1/2$, but with the
opposite orientation induced by complex conjugation. The lemma
follows.
\end{proof}

\section{Related bundles and their homological properties}
\label{sec3}

If $n \geq 2$, the abelianization of $B_n$ is isomorphic to $\mathbb
Z$ with one choice of epimorphism given by $A\co B_n \to\ \mathbb Z$
for $A(\sigma_i) = 1$, $i \geq 1$. The stablization maps $s\co  B_n
\to\ B_{n+1}$ send each $\sigma_i$ to $\sigma_i$ for $ 1 \leq i \leq
n$. These maps are compatible with abelianization via the following
commutative diagram.
\[
\begin{CD}
B_n    @>{s}>>  B_{n+1} @>{s}>>  B_{\infty}   \\
@VV{A}V             @VV{A}V   @VV{A}V         \\
\mathbb Z  @>{\text{Identity}}>>  \mathbb Z @>{\text{Identity}}>>  \mathbb Z
\end{CD}
\]
In addition, there is a map $\theta\co  K(B_{\infty},1) \to\ \Omega^2
S^3$ which induces a homology isomorphism with any simple field
coefficients and with acyclic homotopy fibre denoted $F$
(see Cohen \cite{C,C1}, May \cite{May}, Segal \cite{Segal}). The map
$\theta$ is defined via
\begin{enumerate}
  \item the action of the little $2$--cubes \cite{May} giving a map from the unordered
configuration space $Conf(\mathbb R^2,n)/\Sigma_n \to \Omega^2_{(n)}
S^2$ the path-component of $\Omega^2 S^2$ given by maps having
degree $n$
  \item together with the homotopy equivalence
$$\Omega^2(\eta)\co \Omega^2 S^3 \to \Omega^2_{(0)}S^2$$
where the map $\eta$ denotes the Hopf map $\eta\co S^3 \to S^2$, and
\item passage to colimits.
\end{enumerate} Thus there is a homotopy commutative diagram where $f\co \Omega^2 S^3
\to\ S^1$ induces an isomorphism on the level of fundamental groups.
\[
\begin{CD}
K(B_{\infty},1) @>{\theta}>>  \Omega^2 S^3  \\
@VV{BA}V             @VV{f}V        \\
K(\mathbb Z,1)  @>{\text{Identity}}>>  K(\mathbb Z,1).
\end{CD}
\] This information is recorded in the next lemma.

\begin{lem}
\begin{enumerate}
\item There is a map $K(B_{\infty},1) \to\ \Omega^2 S^3$ which
induces a homology isomorphism, and with acyclic homotopy fibre $F$.

\item The action of $B_n$ on $\mathbb Z[t^{\pm 1}]$ factors
through the abelianization map $A\co B_n \to \mathbb Z$, and thus the
local coefficient system in homology for the fibration
$$EB_n \times_{B_n} Y  \to BB_n$$
factors through the isomorphism
$\pi_1 \Omega^2S^3 \to \pi_1 S^1$.
\end{enumerate}
\end{lem}

Notice that there are several bundles described above over a base
$B$ where
\begin{enumerate}
  \item the base $B$ runs over $BB_n = K(B_{n},1)$, $B_{\infty}= K(B_{\infty},1)$,
  $\Omega^2 S^3$, or $S^1$,
  \item these bundles all have fibre $Y$, and
  \item the action of the fundamental group of $B$ on the homology
  of $Y$ factors through the action of $\pi_1(S^1)$ on $H_*(Y)$.
\end{enumerate}

In what follows below, the notation
$\Gamma(B)$ is used generically to denote any of these bundles.
Thus there is a bundle $\Gamma(S^1)$ over $S^1$ with fibre $Y$
together with a commutative diagram of bundles obtained from pulling
back the bundle over $S^1$.
\[
\begin{CD}
Y  @>{1}>>        Y               @>{1}>>  Y             @>{1}>>  Y    \\
   @VV{i}V          @VV{i}V      @VV{i}V              @VV{i}V \\
\Gamma(BB_n)  @>>>  \Gamma(BB_{\infty})  @>{}>> \Gamma(\Omega^2 S^3) @>{}>> \Gamma(S^1) \\
    @VV{p}V      @VV{p}V      @VV{p}V              @VV{p}V \\
  BB_n   @>>{}>  BB_{\infty}  @>>{}>  \Omega^2 S^3 @>>{}>  S^1
\end{CD}
\]
Recall that there is a homotopy equivalence
$$\Omega^2 S^3 \to\ \Omega^2 {S^3}\<3\> \times  S^1$$
(where $S\<3\>$
denotes the $3$--connected cover of the $3$--sphere, and where this
decomposition is not multiplicative ). Since $\Omega^2 {S^3}\<3\> $ is
simply-connected, the projection $p\co \Omega^2 S^3\<3\> \times S^1 \to
S^1$ induces an isomorphism on the level of fundamental groups.

\begin{thm}\label{thm:more bundles}
There are homotopy equivalences

\begin{enumerate}
  \item $\Omega^2 S^3\<3\> \times \Gamma(S^1) \to \Gamma(\Omega^2
  S^3\<3\> \times S^1)$,
  \item $\Gamma(\Omega^2 S^3) \to  \Gamma(\Omega^2 S^3\<3\> \times S^1)$,
  and
  \item $\Gamma(\Omega^2 S^3) \to  \Omega^2 S^3\<3\> \times \Gamma(S^1)$.
\end{enumerate}

Furthermore, the following properties are satisfied.
\begin{enumerate}
\item[(4)] The map of bundles
$$\Gamma(BB_{\infty}) \to\ \Gamma(\Omega^2 S^3)$$
induces a homology isomorphism with any simple coefficients.

\item[(5)] With simple field coefficients $\mathbb F$, there are
isomorphisms
$$H_*( \Gamma(S^1); \mathbb F) \otimes_{\mathbb F}
  H_*( \Omega^2 S^3 \<3\>; \mathbb F)
  \to\ H_*(\Gamma(\Omega^2 S^3);\mathbb F).$$

\item[(6)] The space $\Gamma(S^1)$ is a $K(\pi,1)$ where $\pi$ is a free
group on $2$ generators.

\item[(7)] If $i \geq 1$, there are isomorphisms $ H_{i+1}(EB_n
\times_{B_n}Y; \mathbb Q) \to\ H_i(B_n; \mathbb Q[t^{\pm 1}]) $ for
$i > 0$.

\end{enumerate}

\end{thm}

Observe that \fullref{thm:bundles} follow at once from \fullref{thm:more bundles} which is proven next.

\begin{proof}
Notice that if $\pi\co A \times B$ denotes the natural projection map
with $F\to E\to B$ a fibre bundle, then the pullback to $A \times B$
is homeomorphic to a product $A \times E$. Thus statements (1)--(3)
follow.

Next, observe that there is a map of fibrations
\[
\begin{CD}
\{*\}               @>>>  Y             @>{1}>>  Y    \\
   @VV{}V          @VV{i}V             @VV{i}V \\
F  @>>>  \Gamma(BB_{\infty})  @>{}>> \Gamma(\Omega^2 S^3)  \\
    @VV{1}V      @VV{p}V                 @VV{p}V \\
F   @>>{}>  BB_{\infty}  @>>{}>  \Omega^2 S^3
\end{CD}
\]
where $F$, an acyclic space, is the homotopy theoretic
fibre of the homology isomorphism $BB_{\infty} \to\ \Omega^2 S^3$.
Since the homology of $F$ is trivial, the Serre spectral sequence
for the fibration $\Gamma(BB_{\infty}) \to\ \Gamma(\Omega^2 S^3)$
collapses, and the induced map is an isomorphism in homology with
any simple coefficients. Part (4) follows.

To work out the homology of $\Gamma(\Omega^2 S^3)$, consider the
following morphism of fibrations
\[
\begin{CD}
\{*\}               @>>>  Y             @>>>  Y    \\
   @VV{}V          @VV{i}V             @VV{i}V \\
\Omega^2 S^3\<3\>  @>>>           \Gamma(\Omega^2 S^3)   @>{}>> \Gamma(S^1)  \\
    @VV{1}V      @VV{p}V                 @VV{p}V \\
\Omega^2 S^3\<3\>   @>>{}>  \Omega^2 S^3  @>>{}>  S^1.
\end{CD}
\]
By statement (4), there is a homotopy equivalence $\Gamma(\Omega^2
S^3) \to  \Omega^2 S^3\<3\> \times \Gamma(S^1)$. Since the natural map
$\Gamma(BB_{\infty}) \to\ \Gamma(\Omega^2 S^3)$ induces a homology
isomorphism with any simple field coefficients $\mathbb F$, there is
an induced map
$$\Gamma(BB_{\infty}) \to\ \Omega^2 S^3\<3\> \times \Gamma(S^1)$$
which induces a homology isomorphism with any simple
field coefficients $\mathbb F$. Part (5) follows.

To finish part (6), the space $\Gamma(S^1)$ will be shown to be
homotopy equivalent to a bouquet of two circles. Notice that
$\Gamma(S^1)$ is given by
$$\mathbb R \times_{\mathbb Z}Y.$$
Let $V$ denote the subspace of $\mathbb C $ given by the union of
circles with centers at $(n,0)$ for integers $n$, and with radii all
equal to $1/2$. Notice that $V$ is a subspace of $Y$, and that the
natural inclusion $V \to\ Y$ is both a homotopy equivalence, and
invariant with respect to the action of $\mathbb Z$ for which a
generator sends $z$ to $1+ \chi(z)$. In addition, this action of
$\mathbb Z$ on $V$ is free, and properly discontinuous. Thus
$E(S^1)$ is homotopy equivalent to the orbit space
$$V/ \mathbb Z.$$
Observe that the orbit space $ V/ \mathbb Z$ is the identification
space of a circle with one pair of antipodal points identified. Thus
$\Gamma(S^1)$ is $K(F_2,1)$ where $F_2$ is a free group with two
generators.

To check statement $7$, notice that if $i \geq 1$, there are
isomorphisms $ H_{i+1}(EB_n \times_{B_n}Y; \mathbb Q) \to\ H_i(B_n;
\mathbb Q[t^{\pm 1}]) $ for $i > 0$ by inspection of the Leray
spectral sequence for the natural projection
$$EB_n \times_{B_n}Y \to BB_n,$$
together with the fact that $H_i(B_n; \mathbb Q) = \{0\}$ for $i > 1$.

The theorem follows.
\end{proof}

\begin{prop}\label{prop:isomorphisms}

There are isomorphisms
$$H_*(\Gamma(BB_{\infty}); \mathbb F)\to H_{*}(B_{\infty}; \mathbb F) \oplus
H_{*-1}(B_{\infty}; \mathbb F[t^{\pm 1}]).$$
\end{prop}

\begin{proof}
Recall the previous commutative diagram
\[
\begin{CD}
\Gamma(BB_{\infty})  @>{}>> \Gamma(\Omega^2 S^3)  \\
    @VV{p}V                 @VV{p}V \\
BB_{\infty}  @>>{}>  \Omega^2 S^3,
\end{CD}
\]
together with the homotopy equivalence
$S^1 \times \Omega^2 S^3\<3\> \to \Omega^2 S^3.$ The proposition
follows from the collapse of the Serre spectral sequence for the
fibrations $p$.
\end{proof}

\begin{cor}\label{cor: another isomorphism}
There are isomorphisms
$$H_*(B_{\infty}; \mathbb F[t^{\pm 1}])\to H_{*}(\Omega^2 S^3\<3\>;\mathbb F).$$
\end{cor}

\section{Remarks}
\label{sec4}

This section lists answers for the cohomology of the braid groups
with various natural choices of local coefficient systems. The
common theme here is that the answers are given by classical number
theoretic constructions.

\begin{enumerate}
    \item The work of De Concini, Procesi, and Salvetti \cite{DPS} gives
that cyclotomic fields arise as the cohomology of the braid groups
with coefficients in the determinant of the Burau representation.
    \item Classical work of Eichler \cite{E}, and Shimura  \cite{S} gives that
$$H^1(SL(2,\mathbb Z); \Sym^k(V_2)\otimes \mathbb R)$$
is the vector space of classical modular forms of weight $k+2$ where $V_2$
is the fundamental representation of $SL(2,\mathbb Z)$, and
$\Sym^k(V_2)$ denotes the $k$--fold symmetric power of $V_2$. The
third braid group $B_3$ is the universal central extension of
$SL(2,\mathbb Z)$. An inspection of the natural spectral sequence
gives similar results for $B_n$, $n = 3,4$. For example $H^*(B_3;
\Sym^k(V_2)\otimes \mathbb R)$ is given by
$$H^*(SL(2,\mathbb Z); \Sym^k(V_2)\otimes \mathbb R) \otimes E[v]$$
where $E[v]$ denotes an exterior algebra with a single generator in
degree $1$.

Natural extensions apply to braid groups with more strands.

\item Work of Furusawa, Tezuka, and Yagita \cite{FTY} gives a
computation of the real cohomology of the classifying space of the
diffeomorphism group of the torus. The answer is given in terms of
classical modular forms.

Related work by Cohen \cite{C2} gives that the real cohomology for the
$n$--stranded braid group of a $k$--punctured torus, $Br_n(S^1 \times
S^1 -\{x_1,x_2, \cdots, x_k\})$ is given in terms of classical
modular forms, while a similar result is obtained after tensoring
with the natural sign representation. Fix a dimension $d$, and
stabilize $H^d(Br_n(S^1 \times S^1 -\{x_1,x_2, \cdots, x_k\};
\mathbb R)$ by letting $k$ go to $\infty$. A direct check of the the
dimension of the resulting cohomology groups is that of certain
spaces of Jacobi forms as calculated by Eichler, and Zagier.

\end{enumerate}

Two problems arise naturally in this context.

\begin{enumerate}
    \item What is the homology of $EB_n \times_{B_n} Y$ ?

\item Consider a symplectic representation of $B_n \to Sp(2g,
\mathbb Z).$ What are the cohomology groups $H^*(B_n;
\Sym^k(V_{2g})\otimes \mathbb R)$ for which the local coefficient
system is given by symmetric powers of the symplectic representation
? It is natural to ask whether these admit interpretations analogous
to those given by the Eichler-Shimura isomorphism.

\end{enumerate}

\bibliographystyle{gtart}
\bibliography{link}

\end{document}